\documentclass[12pt,twoside,reqno]{amsart}
\usepackage{tikz}
\usetikzlibrary{calc,decorations.markings,decorations.pathmorphing,fadings,shadings}
\usepackage{amssymb}
\usepackage{amsmath}
\usepackage{amsthm}
\usepackage{hyperref}
\newtheorem{theorem}{Theorem}[section]
\newtheorem{lemma}[theorem]{Lemma}
\newtheorem{proposition}[theorem]{Proposition}
\newtheorem{corollary}[theorem]{Corollary}

\newtheorem{remark}[theorem]{Remark}

\usepackage[ansinew]{inputenc} 
\usepackage{ifthen}
\usepackage{animate}

\newcounter{angle}
\title[Contributions to a Conjecture of Mueller and Schmidt]{Contributions to a Conjecture of Mueller and Schmidt on Thue inequalities}
\author[N.\ Saradha and Divyum Sharma]{N.\ Saradha and Divyum Sharma}
\address{School of Mathematics\\Tata Institute of Fundamental Research, Homi Bhabha Road, Mumbai\,-\,400\,005, INDIA.}
\email{saradha\symbol{64}math.t{i}fr.res.in, divyum\symbol{64}math.t{i}fr.res.in}
\begin{document}

 \begin{abstract}
Let $F(X,Y)=\sum\limits_{i=0}^sa_iX^{r_i}Y^{r-r_i}\in\mathbb{Z}[X,Y]$ be a form of degree $r=r_s\geq 3$, 
irreducible over $\mathbb{Q}$ and having at most $s+1$ non-zero coefficients. Mueller and Schmidt showed
that the number of solutions of the Thue inequality 
\[
 |F(X,Y)|\leq h
\]
is $\ll s^2h^{2/r}(1+\log h^{1/r})$. They \textit{conjectured} that $s^2$ may be replaced by $s$. Let
\[
 \Psi = \max_{0\leq i\leq s} \max\left( \sum_{w=0}^{i-1}\frac{1}{r_i-r_w},\sum_{w= i+1}^{s}\frac{1}{r_w-r_i}\right).
\]
 Then we show that $s^2$ may be replaced by
 $\max(s\log^3s, se^{\Psi})$. We also show that 
 if $|a_0|=|a_s|$ and $|a_i|\leq |a_0|$ for $1\leq i\leq s-1$, then $s^2$ may be replaced by $s\log^{3/2}s$.
  In particular, this is true if $a_i\in\{-1,1\}$.
\end{abstract}
\maketitle

\section{Introduction}
\subsection{Thue inequalities for forms with few coefficients}
Let $F(X,Y)$ be a form of degree $r\geq 3$ with integer coefficients, irreducible over $\mathbb{Q}$ and having 
at most $s+1$ non-zero coefficients. Write
      \begin{equation}\label{poly}
      F(X,Y)=\sum\limits_{i=0}^sa_iX^{r_i}Y^{r-r_i}
      \end{equation}
with $0=r_0<r_1<\ldots<r_s=r$. Let $D$, $H$ and $M$ denote the discriminant, height and Mahler height of $F(X,1)$ 
respectively. For $h\geq 1$, consider the Thue inequality
     \begin{equation}\label{Thue_ineq}
      |F(X,Y)|\leq h.
      \end{equation}
Let $N_F(h)$ denote the number of integer solutions $(x,y)$ of \eqref{Thue_ineq}.
Schmidt \cite{Sc}
 proved that 
      \begin{equation}\label{Schmidt_result}
       N_F(h)\ll\sqrt{rs}\ h^{2/r}(1+\log h^{1/r}).
      \end{equation}
Prior to this result, Bombieri and Schmidt \cite{BS} considered the equation
      \begin{equation}\label{Beqn}
       |F(X,Y)|=h.
      \end{equation}
They showed that 
      \begin{equation}\label{Bombieri_result}
       N_F^{(1)}(h)\ll r^{1+\omega(h)}
      \end{equation}
where $N_F^{(1)}(h)$ denotes the number of solutions $(x,y)$ of \eqref{Beqn} with $\gcd(x,y)$ $=$ $1$ and
$\omega(h)$ denotes the number of distinct prime divisors of $h$.
Since $s\leq r$, it can be easily seen that
\eqref{Schmidt_result} is better than \eqref{Bombieri_result} if $h$ is given and $r$ is large while
\eqref{Bombieri_result} is better than \eqref{Schmidt_result} if $r$ is given and $h$ is large. The proof of inequality
\eqref{Bombieri_result} is based on the Thue-Siegel principle and $p$-adic arguments. Schmidt's result depended also on the
analysis of location of roots of $F(X,1)$. 
 Bombieri modified a conjecture
of Siegel on the inequality \eqref{Thue_ineq} as
     \[
        N_F(h)\leq C(s,h)
     \]
where $C(s,h)$ depends only on $s$ and $h$. (See Mueller and Schmidt \cite[p. 208]{MS2}). 
This was shown to be true in the case $s=1$
by Hyyr\"{o} \cite{Hy}, Evertse \cite{Ev} and Mueller \cite{Mu}. 
The case $s\geq2$ was considered by Mueller and Schmidt in \cite{MS1}
and \cite{MS2}. They proved that 
         \begin{equation}\label{MS_result}
          N_F(h)\ll s^2C_1(r,h)
         \end{equation}
where $C_1(r,h)=h^{2/r}(1+\log h^{1/r})$. 
 In all
the estimates for $N_F(h)$, the factor $h^{2/r}$ in $C_1(r,h)$ is unavoidable. The 
logarithmic factor in $C_1(r,h)$ was improved by Thunder when $h$ is large, see \cite{T1}
and \cite{T2}. In fact, in \cite{T3}, Thunder proved a general result on the number of solutions of decomposable form inequalities (higher dimensional Thue inequalities), which implies that for Thue inequalities \eqref{Thue_ineq}, we have $$N_F(h)\leq C_2(r)h^{2/r},$$ where $C_2(r)$ is an effectively computable, but inexplicit number depending only on $r$.
Further, it is expected that $N_F(h)$ is close to $A_Fh^{2/r}$, where $A_F$ is the area of the set of $(\xi,\eta)\in\mathbb{R}^2$
with $|F(\xi,\eta)|\leq 1$. In \cite{Ma}, Mahler showed the finiteness of $A_F$.
Later, Mueller \&
Schmidt \cite{MS2} and Bean \cite{Be} gave explicit upper estimates for $A_F$. For instance, for any form $F$ with
non-zero discriminant, Bean showed that 
            \[
               A_F\leq 3B\left(\frac{1}{3},\frac{1}{3}\right)<16
            \]
where $B(m,n)$ is the classical beta function. The bound $3B(\frac{1}{3},\frac{1}{3})$ is
attained for forms which are equivalent under $SL_2(\mathbb{R})$ to $XY(X-Y)$. We refer to Evertse \& Gy\H{o}ry \cite{EGY} and Gy\H{o}ry \cite{G2} for more results on Thue equations and Thue inequalities.\\
\\
It was conjectured in \cite{MS2} that it may be possible to replace
the factor $s^2$ in \eqref{MS_result} by $s$. In Theorem \ref{thm_1} below, we show that when the non-zero terms of $F$ are sufficiently far apart, then $s^2$ can be
improved. In Theorem \ref{line}, we improve \eqref{MS_result} for forms with restricted coefficients.
    \begin{theorem}\label{thm_1}
     Let $F(X,Y)$ be given by \eqref{poly}. Put
           \begin{equation}\label{series}
           \Psi= \max_{0\leq i\leq s} \max\left( \sum_{w=0}^{i-1}\frac{1}{r_i-r_w},\sum_{w= i+1}^{s}\frac{1}{r_w-r_i}\right)
           \end{equation}
     and
    \[
    \Phi = \max(\Psi, 3\log\log s).
    \]
Then we have
           \begin{equation}\label{result_1}
               N_F(h)\ll se^{\Phi}C_1(r,h).
           \end{equation}
(In \eqref{series}, an empty sum is taken to be equal to zero.)
    \end{theorem}
 \begin{remark}\label{r4s}
 Suppose $r\leq 4se^{2\Phi}$. Then \eqref{Schmidt_result} implies \eqref{result_1}. Thus, for proving Theorem \ref{thm_1},
 we may assume that 
            \begin{equation}\label{assumption_r}
                r> 4se^{2\Phi}.
            \end{equation}
Further, we may take $s\gg 1$, as otherwise inequality \eqref{MS_result} is sufficient.
\end{remark}   
\noindent Throughout the paper, $c_1,c_2,\ldots$ denote positive absolute constants. We shall illustrate Theorem \ref{thm_1} with some examples below. 
\begin{remark}\label{rk_13}
Let $0\leq i,w\leq s$, $i\neq w$.
 Since $|r_i-r_w|\geq |i-w|$, it follows that $$\Psi\leq \max_{0\leq i\leq s} \max\left( \sum_{n=1}^i\frac{1}{n},\sum_{n=1}^{s-i}\frac{1}{n}\right)\leq\log s+c_1,$$
 (see \cite[p. $55$]{PNT}). Thus we get \eqref{MS_result}.
\end{remark}
\begin{remark}
 We give some instances when $s^2$ in \eqref{MS_result} can be improved.
      \begin{enumerate}
       \item[(i)] Suppose  $|r_i-r_w|\geq c_2|i-w|$ with $c_2> 1$. Then $\Psi\leq \frac{1}{c_2}\log s+c_3$. Hence
		    \[
			N_F(h)\ll s^{1+\frac{1}{c_2}}\ C_1(r,h).
		    \]
        \item[(ii)]  Suppose  $|r_i-r_w|\geq \frac{1}{3}|i-w|\log|i-w|.$
        Then 
        \[
        \Psi\leq 1+\max_{0\leq i\leq s} \max\left( \sum_{n=2}^i\frac{3}{n\log n}, \sum_{n=2}^{s-i}\frac{3}{n\log n}\right)\leq 3\log \log s+c_4,
        \]
        (see \cite[p. $70$]{PNT}).
        Hence
		    \[
			N_F(h)\ll s \log^3s\ C_1(r,h).
		    \]
      \end{enumerate}
\end{remark}

In a different direction, we impose restrictions on the coefficients of $F$ and improve \eqref{MS_result}.
\begin{theorem}\label{line}
  Suppose that the coefficients of
     $ F(X,Y)$
 satisfy 
\begin{equation}\label{st_line_1}
 \left|\frac{a_0}{a_s}\right|^{1/r_s}\leq \left|\frac{a_0}{a_i}\right|^{1/r_i}\textrm{ for } i=1,\ldots,s-1.
\end{equation}
If $r\geq s\log^3s$, then
\[
 N_F(h)\ll s (\log s)\ h^{2/r}.
\]
\end{theorem}
\begin{remark}\label{line2}
 If $r<s\log^3 s$, then by \eqref{Schmidt_result} we have
 \[
   N_F(h)\ll s\log^{3/2}s \ C_1(r,h).
 \]
\end{remark}
As an immediate consequence of Theorem \ref{line} and Remark \ref{line2} we get
\begin{corollary}
 Suppose that $|a_0|=|a_s|=H$, where $H$ is the height of $F$. Then
 \[
   N_F(h)\ll s\log^{3/2}s \ C_1(r,h).
 \]
 In particular, the above estimate is valid if the coefficients $a_i$ assume only the values $\pm 1$.
\end{corollary}


 Theorem \ref{thm_1} is a consequence of the following result.
	  \begin{proposition}\label{prop}
	     Let $F(X,Y)$ be given by \eqref{poly}. Then
	         \[
	          N_F(h)\ll s\left(\frac{\log s}{\Phi}+e^{\Phi+c_5(\log^3 s)e^{-\Phi}}\right)C_1(r,h).
	         \]
	  \end{proposition}
Let $X_1$ and $X_2$ be positive numbers. Divide the solutions $(x,y)$ of \eqref{Thue_ineq}
into three sets according as
\begin{center}
 $\max(|x|,|y|)>X_1$;  $\max(|x|,|y|)\leq X_1$ and $ \min(|x|,|y|)\geq X_2$;\\
 $ \min(|x|,|y|)< X_2$.
\end{center}
In \cite{MS2}, the solutions in these sets were called large, medium and small respectively.
Denote the number of \textit{primitive} solutions, i.e. solutions $(x,y)$ with $\gcd(x,y)=1$,
in these sets by $P_{\ell ar}(X_1)$, $P_{med}(X_1,X_2)$ and $P_{sma}(X_2)$ respectively.
If $X_2>X_1$, put $P_{med}(X_1,X_2)=0$.
We can bound $N_F(h)$ by finding estimates for these quantities (See Section \ref{last}).
 Mueller and Schmidt had shown that 
\begin{equation}\label{la}
 P_{\ell ar}(Y_W)\ll s
\end{equation}
where $Y_W$ is as given in \eqref{YL} below.
In Theorem \ref{thm_2}(i), we have lowered the value of $Y_W$ and obtained the same conclusion
as in \eqref{la}. For medium solutions, we 
 analyze precisely the large derivatives which play an important role in obtaining good
 rational approximations. (See Lemma \ref{existence_u_v}).
The small solutions are handled as in \cite{MS2}.
\subsection{Large solutions}
It was noted by Brindza, Pint\'{e}r, van der Poorten \& Waldschmidt \cite{BPPW} that it seems likely that almost
all the solutions of \eqref{Beqn} are small, and around $h^{1/r}$ provided that $h$ is large compared to $H$.
In this section, we present some results supporting this observation.
We use some notations from \cite{BS}.
Choose numbers $a,b$ with $0<a<b<1$. Define
       \[
         t=\sqrt{2/(r+a^2)},\  \lambda=2/((1-b)t),
       \]
 \begin{equation}\nonumber
  \delta=\frac{(r+b^2)t^2-2}{r-1},\  A=\frac{1}{a^2}\left(\log M+\frac{r}{2}\right).
\end{equation}
Further, we put
\[
 B=\frac{2^{r}r^{r/2}M^rh}{\sqrt{|D|}},
\]
\[
           Y_E=(2B\sqrt{|D|})^{1/(r-\lambda)}(4e^A)^{\lambda/(r-\lambda)},
          \]
          \begin{equation}\label{YL}
          Y_G=(2B)^{\frac{1}{r-2}+\frac{1}{r^2}}\ \textrm{ and }\  Y_W=R_1^{1/(r-\lambda)}Y_E
         \end{equation}
         where
 \begin{equation}\label{R_1}
        R_1=e^{800\log^3r}.
        \end{equation}
Since $|D|\leq r^rM^{2r-2}$ (\cite[Theorem $1$]{Ma1}), we have $B>1$.
We use this notation without any further mention in the rest of the paper.
It was shown by Bombieri \& Schmidt \cite{BS} that
$P_{\ell ar}(Y_E)\ll r$.
Improving upon \cite{BPPW}, Gy\H{o}ry \cite{G1} showed that $P_{\ell ar}(Y_G)\leq 25r$. (In fact, he proved that $P_{\ell ar}(Y_G)\leq 5r$ if $r$ is sufficiently large.)
Since $Y_G$ is much smaller than $Y_E$, the latter result is better than the former.
Further, Schmidt \cite{Sc} showed that  $$P_{\ell ar}(Y_E)\ll\sqrt{rs}$$ 
and as mentioned already in \eqref{la}, Mueller \& Schmidt obtained that 
     \[
      P_{\ell ar}(R_1^{1/(r-\lambda)}Y_E)\ll s.
     \]
We shall improve the result of Gy\H{o}ry mentioned above in the following theorem. Further,
using a result of Mignotte \cite{Mi} on the distribution of
zeros of polynomials (see Lemma \ref{mig}), we give different upper 
bounds for large primitive solutions.
  \begin{theorem}\label{thm_2}
  Let $F(X,Y)$ be as in \eqref{poly}. Then
  \begin{enumerate}
   \item[(i)] $P_{\ell ar}(Y_G\ R_1^{\frac{1}{r-2}+\frac{1}{r^2}})\ll s$.
   \item[(ii)] Put
  \[
   \Delta=\frac{\sqrt{3|D|}}{2r^{\frac{r+2}{2}}M^{r-1}} \textrm{ and } R_2=1+\frac{Mr}{2\Delta}.
  \]
  Then
  \begin{equation*}\label{thm_2_result}
         P_{\ell ar}(Y_G\ R_2^{\frac{1}{r-2}+\frac{1}{r^2}})\ll \sqrt{r(\log r+\log M)}.
      \end{equation*}
%
        \end{enumerate}
\end{theorem}
\begin{remark}
 Note that $$Y_G\ R_1^{\frac{1}{r-2}+\frac{1}{r^2}}\ll Y_G$$ and 
 $$Y_G\ R_2^{\frac{1}{r-2}+\frac{1}{r^2}}\ll Y_G^2.$$
 When the coefficients $a_i$ of $F$ assume only the values $\pm 1$, then by a result of Borwein and Erdelyi \cite{BE},
 there can be $\ll\sqrt{r}$ roots in the strip $|\textrm{Im}(z)|\leq\frac{2}{\sqrt{r}}$. Using this fact in place of the result
 of Mignotte, it is possible to show that
 \[
   P_{\ell ar}(Y_G\ R_3^{\frac{1}{r-2}+\frac{1}{r^2}})\ll \sqrt{r}
 \]
where
\[
 R_3=1+M\sqrt{r}.
\]
Note also that
\[
  Y_G\ R_3^{\frac{1}{r-2}+\frac{1}{r^2}}\ll Y_G^2.
\]
\end{remark}
\section{Counting elements of a set}
On many occasions, we need to estimate the cardinality of a set whose elements satisfy certain properties which are termed as
gap principle in the theory of Thue equations. In the following lemma, we give two instances in a formal setup.
\begin{lemma}\label{count_lemma}
 Let $n\geq 2$ and let $U=\{u_1,\ldots,u_{n}\}$ be a set together with a map $T:U\rightarrow \mathbb{R}^{*}$
 such that 
 \[
 A_1\leq T(u_1)\leq T(u_2)\leq\ldots\leq T(u_{n})
 \]
and
\begin{equation}\label{count_1}
 T(u_i)\geq\beta \ T(u_{i-1})^{\gamma} \textrm{ for } 2\leq i\leq n,
\end{equation}
where $\beta>0$, $\gamma\geq 2$. Let
\[
 \kappa=\begin{cases}
           2 &\mbox{ if } \beta> 1\\
           1 &\mbox{ if } \beta\leq 1.
          \end{cases}
\]
\begin{enumerate}
 \item[(i)] Suppose that $T(u_n)\leq B_1$ and $A_1\beta^{1/(\kappa(\gamma-1))}>1$. Then
\begin{equation}\label{gap_1}
 n\leq 1+\frac{1}{\log\gamma}\log\left(\frac{\log B_1}{\log A_1+(\log\beta)/(\kappa(\gamma-1))}\right).
\end{equation}
%
\item[(ii)] Suppose that
\begin{equation}\label{count_3}
 T(u_n)\leq (\eta_1\ T(u_1))^{\eta_2}
\end{equation}
with $\eta_1>1$, $1<\eta_2\leq\gamma^{n-1}$. Let $\beta\leq 1$ and $A_1\geq\left(\frac{\eta_1^{\mu}}{\beta}\right)^{1/\nu}$
with $1\leq \mu<\nu<\gamma-1$. Then 
\begin{equation}\label{gap_2}
 n\leq 1+\frac{1}{\log\gamma}\log\left(\eta_2\max\left(\frac{\mu+\nu}{\mu},\frac{1}{1-\nu/(\gamma-1)}\right)\right).
\end{equation}
\end{enumerate}
\end{lemma}
\begin{proof}
From \eqref{count_1}, by induction, we get
 \begin{eqnarray}
\nonumber  T(u_{n})&\geq& \beta^{1+\gamma+\ldots+\gamma^{n-2}}T(u_1)^{\gamma^{n-1}}\\
 &\geq& (\beta^{1/(\kappa(\gamma-1))}T(u_1))^{\gamma^{n-1}}. \label{count_5}
\end{eqnarray}
 $(i)$ Since $T(u_n)\leq B_1$, \eqref{count_5} implies that
 \[
 1< (\beta^{1/(\kappa(\gamma-1))}T(u_1))^{\gamma^{n-1}}\leq B_1.
 \]
Taking logarithm and using $T(u_1)\geq A_1$, we find
\[
 \gamma^{n-1}\leq \frac{\log B_1}{\log A_1+(\log\beta)/(\kappa(\gamma-1))}.
\]
Since $\gamma\geq 2$, $n\geq 2$, the right hand side of the above inequality is $> 1$. Taking logarithm
once again, we get \eqref{gap_1}.\\
$(ii)$ From \eqref{count_3} and \eqref{count_5}, we get
\[
 (\beta^{1/(\gamma-1)}T(u_1))^{\gamma^{n-1}}\leq  (\eta_1\ T(u_1))^{\eta_2}.
\]
Taking logarithms and using $\eta_2\leq\gamma^{n-1}$, $T(u_1)\geq A_1\geq\left(\frac{\eta_1^{\mu}}{\beta}\right)^{1/{\nu}}$,
we get
\begin{eqnarray*}
 \gamma^{n-1}&\leq&\frac{\eta_2\left((1+\frac{\mu}{\nu})\log\eta_1+\frac{1}{\nu}\log\frac{1}{\beta}\right)}{\frac{\mu}{\nu}\log\eta_1+(\frac{1}{\nu}-\frac{1}{\gamma-1})\log\frac{1}{\beta}}\\
 &\leq& \eta_2\max\left(\frac{\mu+\nu}{\mu},\frac{1}{1-{\nu}/(\gamma-1)}\right)
\end{eqnarray*}
since $\frac{\rho_1+\rho_2}{\rho_3+\rho_4}\leq\max(\frac{\rho_1}{\rho_3},\frac{\rho_2}{\rho_4})$ for positive values of
$\rho_1,\rho_2,\rho_3,\rho_4$. Again taking logarithms in the
above inequality, we get \eqref{gap_2}.
%
\end{proof}

\section{A small set of roots close to solutions}
Let $S$ be any finite set of complex numbers and let $\xi$ be a real number.
Define the distance of $\xi$ from $S$, denoted by $d(S,\xi)$, as
    \[
     d(S,\xi)=\min\limits_{\eta\in S}|\xi-\eta|.
    \]
Thus if $\xi\in S$, then $d(S,\xi)=0$. 
In the ensuing discussions, we specialize $S$ as the set of roots
$\alpha_1,\ldots,\alpha_r$ of $$f(Z)=F(Z,1).$$ 
These discussions are valid if we replace $S$ by $S^*$, which is the set of roots $\beta_1,\ldots,\beta_r$ of $F(1,Z)$.
Note that $\{\beta_1,\ldots,\beta_r\}=\{\alpha_1^{-1},\ldots,\alpha_r^{-1}\}$ and $F(1,Z)$ has the same discriminant,
height and Mahler height as $F(Z,1)$.
In \cite{MS2}, the following result was shown.
\begin{lemma}\cite[Lemma 7]{MS2}\label{set_S}\\
 There is a set $S_1\subseteq S$ with $|S_1|\leq 6s+4$ such that for any real $\xi$, we have
    \[
    d(S_1,\xi) \leq R_1\ d(S,\xi)
    \]
where $R_1$ is given by \eqref{R_1}.
\end{lemma}
We now state a lemma on the distribution of roots of a polynomial
due to Mignotte \cite{Mi}.
\begin{lemma}\label{mig}
 Let $P$ be an irreducible polynomial of degree $d$, with integer coefficients. Let $V$ be a
 sector (open or closed) in the complex plane, centred at the origin with central angle $2\pi\theta$,
 where $0\leq\theta\leq 1$. Then the number $N(V)$ of roots of $P$ in $V$ satisfies
 \[
  N(V)\leq 2\theta d+c_6\sqrt{2d(1.5\log(2d)+2\log M(P))},
 \]
 where $M(P)$ is the Mahler height of $P$.
\end{lemma}
We apply Lemma \ref{mig} to get a subset of $S$ having a property similar to the set $S_1$ in
 Lemma \ref{set_S}.
\begin{lemma}\label{Mignotte_distbn}
 There is a set $S_2\subseteq S$ with $|S_2|\ll \sqrt{r(\log r+\log M)}$ such that for any real $\xi$, we have
     \[
    d(S_2,\xi) \leq R_2\ d(S,\xi)
    \]
 where $R_2$ is as given in Theorem \ref{thm_2}(ii).
\end{lemma}
\begin{proof}
 Define
    \[
     \Delta(f)=\min\limits_{1\leq i<j\leq r}|\alpha_i-\alpha_j|.
    \]
By \cite[Theorem 2]{Ma1}, we have
    \[
     \Delta(f)>2\Delta,
    \]
where $\Delta$ is as in Theorem \ref{thm_2}(ii).
    \begin{figure}[h]\label{F:}
    \begin{tikzpicture}
    \draw[thick] (-3,1)--(3,1);
    \shade[top color=white,bottom color=gray] (-3,1) rectangle +(6,2);
    \draw[thick] (-3,-1)--(3,-1);
    \shade[top color=gray,bottom color=white] (-3,-1) rectangle +(6,-2);
    \draw[<->,dotted] (-3,0)--(3,0) node[below] {$x$}; 
    \draw[<->,dotted] (0,-3)--(0,3) node[left] {$y$}; 
    \draw (0,0) circle (2.5cm);
    \node[below] at (2,0) {$V_1$};
    \node[below] at (-2,0) {$V_2$};
     \node[below] at (1.5,-2) {$\Gamma$};

      \draw (2.86,1.25) -- (-2.86,-1.25);
      \draw (-2.86,1.25) -- (2.86,-1.25);
     \draw (1,0) arc (0:15:1.5);
    \node[right] at (0.95,0.3) {\footnotesize$\frac{2 \pi}{r}$};

     \draw (0.95,0) arc (0:-15:1.5);
    \node[right] at (0.9,-0.3) {\footnotesize$\frac{2 \pi}{r}$};

    \node[right] at (-3.7,-1) 
            {\footnotesize$L_2$};
             \node[right] at (-3.7,1) 
            {\footnotesize$L_1$};
    %
    \end{tikzpicture}
    \caption{}
    \end{figure}
    
Therefore there can be at most one root of $f$ inside the circle $\Gamma$ centred
at the origin and having radius $\Delta$. Now consider the sector $V_1$ shown in 
Figure $1$. It is centred at the origin and has central angle $4\pi/r$ . Let $N(V_1)$ denote
the number of zeros of $f$ inside the sector $V_1.$ By Lemma \ref{mig} with $M(P)=M(f)=M$, we
have
    \[
     N(V_1)\leq 4+c_6\sqrt{2r(1.5\log(2r)+2\log M)}\ll \sqrt{r(\log r+\log M)}.
    \]
 We also have $N(V_2)\ll \sqrt{r(\log r+\log M)}$. If $m$ denotes the number
of roots inside $V_1\cup V_2\cup \Gamma$, then $m\ll \sqrt{r(\log r+\log M)}$. The remaining roots either
lie above the line $L_1$ or below the line $L_2$. This means that except for $m$ roots, say $\alpha_1,\ldots,\alpha_m$,
all the other roots $\alpha$ satisfy 
\[
  |\textrm{Im}(\alpha)|\geq\Delta\left|\sin\left(\frac{2\pi}{r}\right)\right|.
\]
 Let $\xi\in\mathbb{R}$.
Suppose
    \[
     d(S,\xi)=|\xi-\alpha_{i_0}|,
    \]
where $i_0\notin\{1,\ldots,m\}$. Then 
     \begin{eqnarray*}
      |\xi-\alpha_1|&\leq& |\xi-\alpha_{i_0}|+|\alpha_{i_0}-\alpha_1|\\
                    &\leq& |\xi-\alpha_{i_0}|\left(1+\frac{2M}{|\xi-\alpha_{i_0}|}\right).
     \end{eqnarray*}
Now $|\xi-\alpha_{i_0}|\geq|\textrm{Im}(\alpha_{i_0})|\geq\Delta|\sin(\frac{2\pi}{r})|$. Thus for $r\gg 1$, we have
\[
    |\xi-\alpha_{1}|\leq |\xi-\alpha_{i_0}|\left(1+\frac{Mr}{2\Delta}\right).
\]
Thus
\[
  |\xi-\alpha_{1}|\leq R_2 |\xi-\alpha_{i_0}|.
\]
Hence there is a set $S_2$ of $m$ roots such that for any real number $\xi$, we have
   \[
     d(S_2,\xi) \leq R_2\ d(S,\xi).
   \]
This completes the proof of the lemma.
\end{proof}
\section{Proof of Theorem \ref{thm_2}}
We begin by stating a lemma of Lewis and Mahler \cite[Lemma $1$]{L} on Diophantine approximation.
\begin{lemma}\label{favour_root}
If $(x,y)$ is a solution of \eqref{Thue_ineq} with $y\neq 0$ and $H(x,y)=\max(|x|,|y|)>B^{1/r}$, then 
\begin{equation}\nonumber
   d\left(S,\ \frac{x}{y}\right) \leq \frac{B}{H(x,y)^r}.
\end{equation}
\end{lemma}
\qed\\
Combining this lemma with Lemma \ref{Mignotte_distbn}, we obtain the following result.
\begin{lemma}\label{less than half}
There is a set $S_2\subseteq S$ with $|S_2|\ll \sqrt{r(\log r+\log M)}$ such that
 \begin{equation}\nonumber
   d\left(S_2,\ \frac{x}{y}\right) \leq \frac{BR_2}{H(x,y)^r}
\end{equation}
for any solution $(x,y)$ of \eqref{Thue_ineq} with $y\neq 0$ and $H(x,y)>B^{1/r}$.
\end{lemma}
\qed\\
Let $a,b,\delta, t,\lambda$ and $A$ be as in Section $1.3$. 
Following \cite{BS}, we say that a rational number $x/y$ is a $\emph{very good approximation}$ to 
 $\alpha\in S $ if
\begin{equation}
\nonumber    \left|\alpha-\frac{x}{y}\right|<(4e^{A}H(x,y))^{-\lambda}.
\end{equation}
We now state the Thue-Siegel principle as given in \cite[p. 74]{BS}.
\begin{lemma}\label{thuesiegelprinciple}
Let $\alpha\in S$. If $x/y,x'/y'$ are two very good approximations to $\alpha$, then
\begin{equation}
\nonumber   \log(4e^{A})+\log H(x',y')\leq \delta^{-1}\{\log(4e^{A})+\log H(x,y)\}.
\end{equation}
\end{lemma}
\qed\\
\textbf{Proof of Theorem \ref{thm_2}}\\
$(i)$ We argue as in \cite{BS} and \cite{G1}. Let $\alpha_i\in S_1$, where $S_1$ is given by Lemma \ref{set_S}. We first count all
large primitive solutions of 
\eqref{Thue_ineq} which are closest to $\alpha_i$.
Let
\begin{equation}
\nonumber I_i=\left\{(x,y)\ : \ y > Y_G R_1^{\frac{1}{r-2}+\frac{1}{r^2}} , \ \gcd(x,y)=1, \left|\alpha_i-\frac{x}{y}\right|=  d\left(S_1,\ \frac{x}{y}\right)  \right\}.
\end{equation}
In Lemma \ref{count_lemma} we take $U=I_i$.  Suppose that $|I_i|=n$.
Enumerate the elements of $I_i$ as $(x_1,y_1),(x_2,y_{2}), \ldots ,(x_n,y_n)$ with 
$H(x_1,y_1)\leq H(x_2,y_{2})\leq \cdots \leq H(x_n,y_n).$ Take $T((x_j,y_j))=H(x_j,y_j)$.
Then
 \begin{equation}
 \nonumber  \frac{1}{y_j y_{j+1}} \leq  \left|\frac{x_{j+1}}{y_{j+1}}-\frac{x_j}{y_j}\right|    \leq  \left|\frac{x_{j+1}}{y_{j+1}}-\alpha_i\right|+\left|\alpha_i-\frac{x_j}{y_j}\right|    \leq  \frac{2BR_1}{H(x_j,y_j)^r}
 \end{equation}
 by Lemmas \ref{set_S} and \ref{favour_root}. Thus,
 \begin{equation}
 \nonumber   H(x_{j+1}, y_{j+1}) \geq H(x_j,y_j)^{r-1}/(2BR_1).
 \end{equation}
Hence \eqref{count_1} is valid with $$\beta=1/(2BR_1)<1,\ \gamma=r-1.$$
We first estimate the number $n_1$ of elements $(x,y)$ in $I_i$ satisfying
$$H(x,y)\geq \left(\frac{R_1}{\sqrt{|D|}}\right)^{\frac{1}{r-\lambda}}Y_E=(2BR_1)^{\frac{1}{r-\lambda}}(4e^A)^{\frac{\lambda}{r-\lambda}}.$$
All such solutions
 are very good approximations to $\alpha_i$.
Hence by Lemma \ref{thuesiegelprinciple}, we have
\[
 \log T((x_n,y_n))\leq \delta^{-1}\ (\log(4e^{A})+\log T((x_1,y_1)))
\]
i.e.
\[
 T((x_n,y_n))\leq (4e^{A}T((x_1,y_1)))^{1/\delta}.
\]
In Lemma \ref{count_lemma}(ii), take 
\[
 A_1=(2BR_1)^{\frac{1}{r-\lambda}}(4e^A)^{\frac{\lambda}{r-\lambda}},
\]
$$\eta_1=4e^A,\ \eta_2=1/\delta, \ \mu=\lambda\ \textrm{ and }\ \nu=r-\lambda.$$
As mentioned in Remark \ref{r4s}, we may take $r$ sufficiently large.
Note that $1/\delta\sim r$. Hence
$1<\eta_2\leq\gamma^{n-1}$ for $n\geq 3$. Also, $\mu=\lambda\sim\sqrt{r}$. Hence $1\leq\mu<\nu=r-\lambda<r-2.$ Thus
\[
n_1\ll 1+\frac{1}{\log(r-1)}\log\left(r\max\left(\frac{r}{\lambda},\frac{r-2}{\lambda-2}\right)\right)\ll 1.
\]
Next, we
estimate the number $n_2$ of solutions $(x,y)$ satisfying
\[
Y_G R_1^{\frac{1}{r-2}+\frac{1}{r^2}} < H(x,y)\leq \left(\frac{R_1}{\sqrt{|D|}}\right)^{1/(r-\lambda)}Y_E\leq (2BR_1)^{\frac{1+\lambda/a^2}{r-\lambda}}.
\]
Applying Lemma \ref{count_lemma}(i) with $$A_1=Y_G R_1^{\frac{1}{r-2}+\frac{1}{r^2}}=(2BR_1)^{\frac{1}{r-2}+\frac{1}{r^2}},$$
$$B_1=(2BR_1)^{\frac{1+\lambda/a^2}{r-\lambda}}, \ \beta=1/(2BR_1), \ \kappa=1\ \textrm{ and }\ \gamma=r-1,$$ we obtain that $n_2\ll 1$. Thus
\[
 |I_i|=n_1+n_2\ll 1.
\]
As this inequality is true for every root $\alpha_i\in S_1$ and $|S_1|\ll s$, we obtain that the number of
solutions $(x,y)$ with $H(x,y)> Y_G R_1^{\frac{1}{r-2}+\frac{1}{r^2}} $ is $\ll s$. 

 The proof of $(ii)$ is similar to the above proof with $R_1$, $S_1$ 
 replaced by $R_2$, $S_2$.
\qed\\
\\
In Sections $5$-$7$, we improve the mechanism for dealing with medium solutions developed by Mueller and Schmidt.
\section{Archimedean Newton polygon and Large derivatives}
Let $F$ be as in \eqref{poly}.
The Archimedean Newton polygon of $f(z)=F(z,1)$ is the lower boundary of the convex hull of the points
      \[
         P_i=(r_i,-\log |a_i|),\ i=0,\ldots,s.
      \]
We label the vertices of the Newton polygon as 
      \[
          P_0=P_{i(0)}, P_{i(1)},\ldots,P_{i(\ell)}=P_s,
      \]
where $0=i(0)<i(1)<\ldots<i(\ell)=s$. For $k=1,\ldots,\ell$, define
$\sigma(i(k))$ as the slope of the line segment $P_{i(k-1)},P_{i(k)}$.
For $k=0,\ldots,\ell-1$, define
$\sigma^+(i(k))$ as the slope of the line segment $P_{i(k)},P_{i(k+1)}$.
Let $\alpha$ be a root of $f$. We define $K(\alpha)$, $k(\alpha)$ as follows.
If $\sigma^+(i(\ell-1))=\sigma(s)<\log|\alpha|+\Psi+\log 3 $, then put
$K(\alpha)=\ell$. If not, then define $K(\alpha)$ as the least integer $K$
in $0\leq K\leq \ell-1$ with $\sigma^+(i(K))\geq\log|\alpha|+\Psi+\log 3$.
If $\sigma(i(1))=\sigma^+(0)>\log|\alpha|-\Psi-\log 3 $, then put
$k(\alpha)=0$. Otherwise, define $k(\alpha)$ as the largest integer $k$
in $1\leq k\leq \ell$ with $\sigma(i(k))\leq\log|\alpha|-\Psi-\log 3$.
Clearly, $k(\alpha)\leq K(\alpha)$.  In \cite{MS2}, $k(\alpha)$ and $K(\alpha)$
are defined with $\Psi$ replaced by $\log s$. As in \cite[Equation $(6.3)$]{MS2}, we have
\[
            k(\alpha)<K(\alpha)
           \]
for every root $\alpha$ of $f$.

    \begin{lemma} \label{rk_st}                                         
     The coefficients of $F$ satisfy condition \eqref{st_line_1} if and only if the Newton polygon of $F$ is a 
     straight line joining $(0,-\log|a_0|)$ and $(r,-\log|a_s|)$. Further, in this case the height $H$ of $F$ is either
     $|a_0|$ or $|a_s|$.
    \end{lemma}
    \begin{proof}
     The Newton polygon is a straight line if and only if the slope of the line joining 
     $(0,-\log|a_0|)$ and $(r,-\log|a_s|)$ does not exceed the slope of the line joining
     $(0,-\log|a_0|)$ and $(r_i,-\log|a_i|)$ for $i=1,\ldots,s-1$, i.e. for every $i$ with $1\leq i\leq s-1$, we have
     \[
      \frac{-\log|a_s|+\log|a_0|}{r}\leq \frac{-\log|a_i|+\log|a_0|}{r_i}
     \]
     or,
     \[
      \left|\frac{a_0}{a_s}\right|^{1/r}\leq\left|\frac{a_0}{a_i}\right|^{1/r_i}
     \]
    which is condition \eqref{st_line_1}. This proves the first assertion.\\
%
\\
      Condition \eqref{st_line_1} means that for $i=1,\ldots,s-1$, we have
      \[
       |a_s||a_0|^{\frac{r}{r_i}-1}\geq |a_i|^{\frac{r}{r_i}}.
      \]
      If $|a_0|\leq|a_s|$, then the above inequality implies that
      \[
       |a_s|^{\frac{r}{r_i}}\geq |a_i|^{\frac{r}{r_i}}  \textrm{ for }i=1,\ldots,s-1.
      \]
      Thus
      \[
       |a_s|\geq|a_i| \textrm{ for } i=0,\ldots,s-1.
      \]
      Hence $H=|a_s|$.
      Similarly, if $|a_s|\leq|a_0|$, then we get
      \[
       |a_0|^{\frac{r}{r_i}}\geq |a_i|^{\frac{r}{r_i}} \textrm{ for }i=1,\ldots,s-1.
      \]
      Thus $H=|a_0|$. This proves the second assertion.
     \end{proof}
\begin{remark}\label{rk_lemma}
 By Lemma \ref{rk_st}, when the coefficients of $F$ satisfy \eqref{st_line_1}, we have $\ell=1$. 
 Since $k(\alpha)<K(\alpha)$, we find 
 $$k(\alpha)=0 \textrm{ and } K(\alpha)=1$$
 for any $\alpha\in S$.
\end{remark}

As a consequence of Lemma \ref{rk_st}, we get the following result which is a special case of Lemma $1(i)$ of \cite{MS2}.
\begin{lemma}\label{bds}
 Suppose that the coefficients of $F$ satisfy \eqref{st_line_1}. Then for every root $\alpha$ of $f$, we have
 \[
  \frac{1}{2}e^{\sigma}<|\alpha|<2e^{\sigma},
 \]
 where $\sigma$ denotes the slope of the line joining $(0,-\log|a_0|)$ and $(r,-\log|a_s|).$ 
\end{lemma}
\begin{proof}
 By Lemma \ref{rk_st}, the coefficients of $F$ satisfy \eqref{st_line_1} if and only if 
 \[
  \sigma\leq \frac{-\log|a_i|+\log|a_0|}{r_i} \textrm{ for }i=1,\ldots,s-1.
 \]
This implies that 
\[
 |a_s|e^{\sigma r}=|a_0|\geq |a_i|e^{\sigma r_i} \textrm{ for }i=0,\ldots,s-1.
\]
When $z=e^{\sigma}w$ with $|w|\geq 2$, we have
\begin{eqnarray*}
 |f(z)|&=&|a_se^{\sigma r}w^r+a_{s-1}e^{\sigma r_{s-1}}w^{r_{s-1}}+\cdots+a_0|\\
 &\geq& |a_se^{\sigma r}||w^r|-|a_{s-1}e^{\sigma r_{s-1}}||w^{r_{s-1}}|-\cdots-|a_0|\\
 &\geq& |a_se^{\sigma r}|(|w^r|-|w^{r_{s-1}}|-\cdots-1)>0.
\end{eqnarray*}
Therefore if $\alpha$ is a root of $f$, then 
\[
 |\alpha|<2e^{\sigma}.
\]
Now consider the reciprocal polynomial $\hat{f}$ of $f$, i.e. 
\[
 \hat{f}(z)=z^rf(1/z).
\]
The Newton polygon of $\hat{f}$ is the single line with slope $-\sigma$.
Hence every root $\hat{\alpha}$ of $\hat{f}$ satisfies
\[
 |\hat{\alpha}|<2e^{-\sigma}.
\]
Since the roots of $\hat{f}$ are the reciprocals of the roots of $f$, we obtain
\[
 |\alpha|>\frac{1}{2}e^{\sigma}
\]
for every root
$\alpha$ of $f$. This completes the proof of the lemma.

\end{proof}
We shall now use the Newton polygon to prove that for each root $\alpha$ of $f$, there exists $u$ with $1\leq u\leq r$ such
that $|f^{(u)}(\alpha)|$ is large. This will enable us to obtain good rational approximations to $\alpha$ from the solutions of \eqref{Thue_ineq}
(See Lemma \ref{dvt1}).
We introduce some notation. Let $e,h$ be two non-negative integers. Let $(e)_h$ be the Pochhammer symbol defined as
  \[
   (e)_h=\begin{cases}
          0 &\mbox{if } e=0\\
          1 &\mbox{if } h=0\\
          e(e-1)\cdots(e-h+1) &\mbox{otherwise }.
         \end{cases}
  \]
%
For a positive integer $t$, define
   \[
	\Delta_t^-(e)=\begin{pmatrix} 
	(e)_0 \\
	\vdots \\
	(e)_t
	\end{pmatrix}. 
   \]
Further, for $0\leq u\leq t$, let
    \[
	\Delta_{t,u}^-(e)=\begin{pmatrix} 
	(e)_0 \\
	\vdots \\
	(e)_{u-1}\\
	(e)_{u+1}\\
	\vdots\\
	(e)_{t}
	\end{pmatrix}. 
   \]
If $\{a_1,\ldots,a_{t+1}\}$ is a set of positive integers, then
   \begin{equation}\label{vandermonde}
       \det\left(\Delta_t^-(a_1),\ldots,\Delta_t^-(a_{t+1})\right)=\prod\limits_{1\leq i<j\leq t+1}(a_j-a_i).
   \end{equation}
Let $\{b_1,\ldots,b_t\}$ be any set of positive integers. Put
   \[
      E_u^{(t)}=E_u^{(t)}(b_1,\ldots,b_t)=(-1)^{t+u}\det\left(\Delta_{t,u}^-(b_1),\ldots,\Delta_{t,u}^-(b_t)\right).
   \]
Then for any positive integer $e$, we have
   \[
      \sum\limits_{u=0}^t (e)_u E_u^{(t)}(b_1,\ldots,b_t)=\det\left(\Delta_t^-(b_1),\ldots,\Delta_t^-(b_{t}),\Delta_t^-(e)\right).
   \]
We denote $\det\left(\Delta_t^-(b_1),\ldots,\Delta_t^-(b_{t}),\Delta_t^-(e)\right)$ by $D(b_1,\ldots,b_t,e)$. Note that
$D(b_1,\ldots,b_t,e)$ $=0$ whenever $e=b_i$ for any $i$ with $1\leq i\leq t$.
   \begin{lemma}\label{3.1}
     Let $P(z)=p_1z^{e_1}+\ldots+p_mz^{e_m}$, $e_1<\ldots<e_m$, be a polynomial and let 
     $\{b_1,\ldots,b_t\}$ be any set of positive integers. Then
        \[
           \sum_{u=0}^tE_u^{(t)}z^uP^{(u)}(z)=\sum_{i=1}^mp_iz^{e_i}D(b_1,\ldots,b_t,e_i),
        \]
     where $ E_u^{(t)}=E_u^{(t)}(b_1,\ldots,b_t)$.
   \end{lemma}
   \begin{proof}
     Observe that for $0\leq u\leq t$, we have
	    \begin{eqnarray*}
		z^uP^{(u)}(z)&=&(e_1)_u\ p_1z^{e_1}+\ldots+(e_m)_u\ p_mz^{e_m}.
	    \end{eqnarray*}
     Then
            \begin{eqnarray*}
               \sum_{u=0}^tE_u^{(t)}z^uP^{(u)}(z)&=&
              p_1z^{e_1}\sum_{u=0}^t (e_1)_u E_u^{(t)}+\ldots+p_mz^{e_m}\sum_{u=0}^t (e_m)_uE_u^{(t)}\\
              &=&p_1z^{e_1}D(b_1,\ldots,b_t,e_1)+\ldots+p_mz^{e_m}D(b_1,\ldots,b_t,e_m).
            \end{eqnarray*}
    This proves the lemma.    
    \end{proof}
    \begin{corollary}
     Let $\alpha$ be a root of $f(z)$. Then
	\begin{equation}\label{E_i(K)}
	  \sum_{u=1}^{i(K)}E_u^{(i(K))}\alpha^uf^{(u)}(\alpha)=a_{i(K)}D_{i(K)}^{(1)}\alpha^{r_{i(K)}}+\sum_{j>i(K)}a_jD_j^{(1)}\alpha^{r_j},
	\end{equation}
     where $K=K(\alpha)$, $E_u^{(i(K))}=E_u^{(i(K))}(r_0,\ldots,r_{i(K)-1})$, \\ $D_j^{(1)}$ $=$ $D(r_0,\ldots,r_{i(K)-1},r_j)$ and the sum with $j>i(K)$ is taken as zero
     if $K=\ell$.
     Also
        \begin{equation}\label{E_i(k)}
                 \sum_{u=1}^{s-i(k)}E_u^{(s-i(k))}\alpha^uf^{(u)}(\alpha)=a_{i(k)}D_{i(k)}^{(2)}\alpha^{r_{i(k)}}+\sum_{j<i(k)}a_jD_j^{(2)}\alpha^{r_j}
        \end{equation}
     where $k=k(\alpha)$, $E_u^{(s-i(k))}=E_u^{(s-i(k))}(r_{i(k)+1},\ldots,r_{s})$, \\ $ D_j^{(2)}=D(r_{i(k)+1},\ldots,r_{s},r_j)$ and the sum with $j<i(k)$ is taken as zero if $k=0$.
    \end{corollary}
    \begin{proof}
       In Lemma \ref{3.1}, take $P(z)=f(z)$ with $(e_1,\ldots,e_m)=(r_0,\ldots,r_s)$ and
       $(p_1,\ldots,p_m)=(a_0,\ldots,a_s)$. Now \eqref{E_i(K)} follows from the lemma by taking $z=\alpha$, $t=i(K)$,
       $(b_1,\ldots,b_t)=(r_0,\ldots,r_{i(K)-1})$ and using the fact that $D(r_0,\ldots,r_{i(K)-1},r_j)=0$ for $j<i(K)$.
       Similarly, we obtain \eqref{E_i(k)} by taking $z=\alpha$, $t=s-i(k)$,
       $(b_1,\ldots,b_t)=(r_{i(k)+1},\ldots,r_{s})$ and using the fact that $D(r_{i(k)+1},\ldots,r_{s},r_j)=0$ for $j>i(k)$.
    \end{proof}
Note that the sums on the left hand sides of \eqref{E_i(K)} and \eqref{E_i(k)} are non-empty since $0\leq i(k)<i(K)\leq s$.
   \begin{lemma}\label{existence_u_v}
        Let $(x,y)$ be a solution of \eqref{Thue_ineq} with $y\neq 0$. Let $\alpha$ be a root
        of $f$. 
        Then
        \begin{enumerate}
         \item[(i)] There exists $u$ with $1\leq u\leq i(K)$ such that 
                \begin{equation*}
                  |f^{(u)}(\alpha)|\geq \frac{1}{4s}(2s^2r)^{1-s}|a_{i(K)}||\alpha|^{r_{i(K)}-u}.
                \end{equation*}
        \item[(ii)] There exists $v$ with $1\leq v\leq s-i(k)$ such that 
                \begin{equation*}
                  |f^{(v)}(\alpha)|\geq \frac{1}{4s}(2s^2r)^{1-s}|a_{i(k)}||\alpha|^{r_{i(k)}-v}.
                \end{equation*} 
        \end{enumerate}
      \end{lemma}
   \begin{proof}
  Suppose that $K<\ell$. Let $j>i(K)$. Put
      \[
         W_j^{(1)}=\left|\frac{a_{j}D_{j}^{(1)}\alpha^{r_{j}}}{a_{i(K)}D_{i(K)}^{(1)}\alpha^{r_{i(K)}}}\right|.
      \]
   By \eqref{E_i(K)}, we have
   \begin{align}
      \sum\limits_{u=1}^{i(K)}E_u^{(i(K))}\alpha^uf^{(u)}(\alpha)
                 &=a_{i(K)}D_{i(K)}^{(1)}\alpha^{r_{i(K)}}
                  \left(1 +\sum_{j>i(K)} \frac{a_{j}D_{j}^{(1)}\alpha^{r_{j}}}{a_{i(K)}D_{i(K)}^{(1)}\alpha^{r_{i(K)}}}\right). \label{determinant1}
   \end{align}
	From the definition of $D_j^{(1)}$ and  \eqref{vandermonde} , we obtain
	    \begin{eqnarray*}
	      \left|\frac{D_j^{(1)}}{D_{i(K)}^{(1)}}\right|&=&\prod_{w<i(K)}\left(\frac{r_j-r_w}{r_{i(K)}-r_w}\right)
	      =\prod_{w<i(K)}\left(1+\left(\frac{r_j-r_{i(K)}}{r_{i(K)}-r_w}\right)\right).
	    \end{eqnarray*}
	This implies 
	  \begin{eqnarray*}
	    \log\left|\frac{D_j^{(1)}}{D_{i(K)}^{(1)}}\right|&\leq&(r_j-r_{i(K)})\sum_{w<i(K)}\frac{1}{r_{i(K)}-r_w}\\
	    &\leq&(r_j-r_{i(K)})\Psi.
	  \end{eqnarray*}
	Then
	  \begin{eqnarray*}
	    \log W_j^{(1)}&\leq&(r_j-r_{i(K)})\left(\log|\alpha|+\Psi\right)+\log|a_j|-\log|a_{i(K)}|\\
			  &=&(r_j-r_{i(K)})\left(\log|\alpha|+\Psi-\sigma(j,i(K))\right),
	  \end{eqnarray*}
	where $\sigma(j,i(K))$ denotes the slope of the line segment joining $P_{i(K)}$ and $P_j$. 
	From the convexity of the Newton polygon and the definition of $K=K(\alpha)$, we get
	  \[
	    \sigma(j,i(K))\geq\sigma^+(i(K))\geq\log|\alpha|+\Psi+\log 3.
	  \]
	Hence
	  \[
	    \log W_j^{(1)}\leq-(r_j-r_{i(K)}) \log 3
	  \]
	implying that
	  \[
	      W_j^{(1)}\leq3^{-(r_j-r_{i(K)})}.
	  \]
	Thus
	  \[
	      \sum_{j>i(K)} W_j^{(1)}\leq 3^{-1}+3^{-2}+\ldots=\frac{1}{2}.
	  \]
	Using this in \eqref{determinant1}, we obtain
	  \begin{equation}\label{E1}
	    \left|\sum\limits_{u=1}^{i(K)}E_u^{(i(K))}\alpha^uf^{(u)}(\alpha)\right|\geq \frac{1}{2}|a_{i(K)}D_{i(K)}^{(1)}\alpha^{r_{i(K)}}|.
	  \end{equation}
	  It is easy to see that the above inequality holds also for $K=\ell$
         since then the right hand side of \eqref{determinant1} reduces to $|a_{s}D_{s}^{(1)}\alpha^{r}|$.
         Suppose that $k>0$. As above, it can be shown that for $j<i(k)$,
          \[
	    \log\left|\frac{D_j^{(2)}}{D_{i(k)}^{(2)}}\right| \leq(r_{i(k)}-r_j)\Psi
	  \]
	  and hence
	  \begin{equation}\label{E2}
	    \left|\sum\limits_{u=1}^{s-i(k)}E_u^{(s-i(k))}\alpha^uf^{(u)}(\alpha)\right|
		  \geq \frac{1}{2}|a_{i(k)}D_{i(k)}^{(2)}\alpha^{r_{i(k)}}|.
	  \end{equation}
	  This inequality is also true for $k=0$.
	  From \cite[Eqns $(6.12)$ \& $(6.13)$]{MS2} it follows that
		\begin{equation}\label{det1}
		    |E_u^{(i(K))}|\leq 2^s(s^2r)^{s-1}|D_{i(K)}^{(1)}| \textrm{ for } 1\leq u\leq i(K)
		\end{equation}
	    and
	      \begin{equation}\label{det2}
		    |E_u^{(s-i(k))}|\leq 2^s(s^2r)^{s-1}|D_{i(k)}^{(2)}| \textrm{ for } 1\leq u\leq s-i(k).
		\end{equation}
         Substituting \eqref{det1} in \eqref{E1}, we find that
             \[
              |f^{(u)}(\alpha)|\geq \frac{1}{4s}(2s^2r)^{1-s}|a_{i(K)}||\alpha|^{r_{i(K)}-u}
             \]
         for some $u$ with  $1\leq u\leq i(K)$. 
         In a similar manner, \eqref{E2} and \eqref{det2} yield 
         the second part of the lemma.
      \end{proof}
\section{Good rational approximations}
The results in this section correspond to \cite[Section $14$]{MS2}, with minor changes. We include the details for the convenience
of the reader.
Throughout this section, assume that $(x,y)$ is a solution of \eqref{Thue_ineq} with $y\neq 0$ and that $\alpha$ is a root
        of $f$ with \[
                      d\left(S,\frac{x}{y}\right)=\left|\alpha-\frac{x}{y}\right|.
                    \]
Let $q$ be the smallest integer with
   \[
     |a_q|=H=\max_{0\leq j\leq s}|a_j|.
   \]
Note that $(q,-\log|a_q|)$ is a vertex of the Newton polygon of $f$.
The following lemma from \cite{MS2} gives a rational approximation to $\alpha$ in terms of the derivatives of $f$.
   \begin{lemma}\cite[Lemma 10]{MS2}\label{dvt1}\\
If $f^{(u)}(\alpha)\neq 0$ for some $u$ in $1\leq u\leq r$, then
   \begin{equation*}\label{dvt}
       d\left(S,\frac{x}{y}\right)\leq\frac{r}{2}\left(\frac{2^rh}{|f^{(u)}(\alpha)y^r|}\right)^{1/u}.
   \end{equation*}
\end{lemma}
 \qed\\
 The following two results are a consequence of the above lemma.
   \begin{lemma}\label{med}
     If $q<i(K)$, where $K=K(\alpha)$, then
       \begin{equation}\label{med2}
     d\left(S,\frac{x}{y}\right)   \leq\frac{1}{H^{(1/u)-(1/r)}}\left(\frac{(rs)^{2s}(6e^{\Psi})^{r}h}{|y|^r}\right)^{1/u},
       \end{equation}
     where $u$ is chosen according to Lemma \ref{existence_u_v}(i).
   \end{lemma}
   \begin{proof} Combining Lemmas \ref{existence_u_v}(i) and \ref{dvt1}, we obtain
       \begin{equation}\label{app}
        \left|\alpha-\frac{x}{y}\right|\leq\left(\frac{2^r(rs)^{2s}h}{|a_{i(K)}||\alpha|^{r_{i(K)}-u}|y|^r}\right)^{1/u}.
       \end{equation}
   Denote $|a_{i(K)}||\alpha|^{r_{i(K)}-u} $ by $\Delta(\alpha,u) $. Thus
        \[
         \log\Delta(\alpha,u)=(r_{i(K)}-u)\log|\alpha|+\log|a_{i(K)}|.
        \]
  From the definition of $K(\alpha)$ it follows that
       \[
        \log|\alpha|\geq\sigma(i(K))-\Psi-\log 3.
       \]
  Since $q<i(K)$, we have $\sigma(i(K))\geq \sigma(q,i(K))$ as the Newton polygon is convex. Thus
       \begin{eqnarray*}
         \log\Delta(\alpha,u)&\geq&(r_{i(K)}-u)\left(\sigma(i(K))-\Psi-\log 3\right)+\log|a_{i(K)}|\\
            &\geq&(r_{i(K)}-u)\sigma(q,i(K))+\log|a_{i(K)}|-r\log(3e^{\Psi})\\
            &=&(r_{q}-u)\sigma(q,i(K))+\log|a_{q}|-r\log(3e^{\Psi}).
        \end{eqnarray*}
Since $P_q$ is one of the lowest vertices of the Newton polygon and $q<i(K)$, the slope
 $\sigma(q,i(K))$ is non-negative. Therefore if $r_q\geq u$, we have
      \[
       \log\Delta(\alpha,u)\geq \log|a_{q}|-r\log(3e^{\Psi})=\log H-r\log(3e^{\Psi}).
      \]
If $r_q<u$, we have $\sigma(q,i(K))\leq\sigma(q,s)$. Hence
     \begin{eqnarray*}
      \log\Delta(\alpha,u)&\geq& (r_{q}-u)\sigma(q,s)+\log|a_{q}|-r\log(3e^{\Psi})\\
      &=&\left(1-\frac{u-r_q}{r_s-r_q}\right)\log|a_q|+\frac{u-r_q}{r_s-r_q}\log|a_s|-r\log(3e^{\Psi})\\
      &\geq&\left(1-\frac{u-r_q}{r-r_q}\right)\log|a_q|-r\log(3e^{\Psi})\\
      &\geq&\left(1-\frac{u}{r}\right)\log H-r\log(3e^{\Psi}).
     \end{eqnarray*}
Substituting this in \eqref{app}, we obtain the assertion of the lemma.
   \end{proof}
   \begin{lemma}\label{v2}
     If $x\neq 0$, $|y|\geq 2(rs)^{2s/r}h^{1/r}$ and $i(k)<q$, where $k=k(\alpha)$, then
         \begin{equation}\label{v22}
        d\left(S^{*},\frac{y}{x}\right)\leq\frac{1}{H^{(1/v)-(1/r)}}\left(\frac{(rs)^{2s}(12 e^{\Psi})^{r}h}{|x|^r}\right)^{1/v},
       \end{equation}
      where $v$ is chosen according to Lemma \ref{existence_u_v}(ii) and $S^{*}=\{\alpha^{-1}|\alpha\in S\}$.
   \end{lemma}
   \begin{proof}
    Combining Lemmas \ref{existence_u_v}(ii) and \ref{dvt1}, we obtain
       \begin{equation*}\label{app*}
        \left|\alpha-\frac{x}{y}\right|\leq\left(\frac{2^r(rs)^{2s}h}{|a_{i(k)}||\alpha|^{r_{i(k)}-v}|y|^r}\right)^{1/v}.
       \end{equation*}
   Denote $|a_{i(k)}||\alpha|^{r_{i(k)}-v} $ by $\Delta^*(\alpha,v) $. Then
        \[
         \log(|\alpha|^v\Delta^*(\alpha,v))=r_{i(k)}\log|\alpha|+\log|a_{i(k)}|.
        \]
        If $k=0$, we have $r_{i(0)}=0$ and hence $\log(|\alpha|^v\Delta^*(\alpha,v))\geq 0$. Now suppose
        that $k>0$. From the definition of $k(\alpha)$ it follows that
       \[
        \log|\alpha|\geq\sigma(i(k))+\Psi+\log 3.
       \]
 This implies that
  \begin{eqnarray*}
      \log(|\alpha|^v\Delta^*(\alpha,v))&\geq& r_{i(k)}\sigma(i(k))+\log|a_{i(k)}|\\
      &\geq& r_{i(k)}\sigma(0,i(k))+\log|a_{i(k)}|=\log|a_0|\geq 0.
     \end{eqnarray*}
     Therefore
     \begin{equation*}
        \left|\alpha-\frac{x}{y}\right|\leq|\alpha|\left(\frac{2^r(rs)^{2s}h}{|y|^r}\right)^{1/v}\leq|\alpha|
       \end{equation*}
 by the assumption on $y$. Hence
       \[
        |x|\leq|2\alpha y|.
       \]
Using this we obtain 
\begin{eqnarray*}
        \left|\alpha^{-1}-\frac{y}{x}\right|&=&\left|\frac{y}{x\alpha}\right|\left|\alpha-\frac{x}{y}\right|
        \leq\left(\frac{2^{r}(rs)^{2s}h}{|a_{i(k)}||\alpha|^{r_{i(k)}}}\right)^{1/v}\frac{1}{|x||y|^{(r/v)-1}}\\
        &\leq&\left(\frac{4^r(rs)^{2s}h}{\Gamma(\alpha,v)|x|^r}\right)^{1/v},
       \end{eqnarray*}
where $\Gamma(\alpha,v)=|a_{i(k)}||\alpha|^{-(r-r_{i(k)}-v)}$. Note that $r-r_{i(k)}-v\geq 0$.
       It is enough to show that
       \[
        \Gamma(\alpha,v)\geq (3e^{\Psi})^{-r}H^{1-(v/r)}.
       \]
Since $i(k)<q$, we have $\sigma(i(k),q)\geq\sigma^+(i(k))>\log|\alpha|-\Psi-\log 3$. Thus
    \begin{eqnarray*}
     \log\Gamma(\alpha,v)&\geq&(r-r_{i(k)}-v)(-\sigma^+(i(k))-\log(3e^{\Psi}))+\log |a_{i(k)}|\\
     &\geq&-(r-r_{i(k)}-v)\sigma(i(k),q)+\log |a_{i(k)}|-r\log(3e^{\Psi})\\
     &=&-(r-r_{q}-v)\sigma(i(k),q)+\log |a_{q}|-r\log(3e^{\Psi}).
    \end{eqnarray*}
When $v\leq r-r_q$, we have
     \[
      \log\Gamma(\alpha,v)\geq\log |a_{q}|-r\log(3e^{\Psi})=\log H-r\log(3e^{\Psi})
     \]
as $\sigma(i(k),q)\leq 0$. When $v >r-r_q$, we get
     \begin{eqnarray*}
      \log\Gamma(\alpha,v)&\geq&-(r-r_{q}-v)\sigma(0,q)+\log |a_{q}|-r\log(3e^{\Psi})\\
      &=&(r-r_{q}-v)((\log|a_q|-\log|a_0|)/r_q)+\log |a_{q}|-r\log(3e^{\Psi})\\
      &\geq&((r-v)/r_q)\log|a_q|-r\log(3e^{\Psi})\\
      &\geq&(1-(v/r))\log H-r\log(3e^{\Psi}).
     \end{eqnarray*}
This proves the claim and hence the assertion of the lemma.
   \end{proof}

We combine Lemmas \ref{med}, \ref{v2} and \ref{set_S} to get the following result, which is
analogous to \cite[Lemma 17]{MS2}.
   \begin{lemma}\label{app_lemma}
    There is a set $S_1$ of roots of $F(x,1)$ and a set $S_1^*$ of roots of $F(1, y)$, both
    with cardinalities $\leq 6s+4$, such that any solution $(x,y)$ of \eqref{Thue_ineq} with
    \[ \min(|x|,|y|) \geq 12e^{\Psi}(rs)^{2s/r}h^{1/r} \] either has
	\begin{equation}\label{medium_app}
	    \left|\alpha-\frac{x}{y}\right|\leq\frac{R_1}{H^{\frac{1}{s}-\frac{1}{r}}}\left(\frac{(rs)^{2s}(12e^{\Psi})^rh}{|y|^r}\right)^{1/s}
	\end{equation}
    for some $\alpha\in S_1$ or has
	\begin{equation}\label{29dash}
	    \left|\alpha^*-\frac{y}{x}\right|\leq\frac{R_1}{H^{\frac{1}{s}-\frac{1}{r}}}\left(\frac{(rs)^{2s}(12e^{\Psi})^rh}{|x|^r}\right)^{1/s}
	\end{equation}
	for some $\alpha^*\in S_1^*$.
    \end{lemma}
    \begin{proof}
     Since $\min(|x|,|y|) \geq 12e^{\Psi}(rs)^{2s/r}h^{1/r}$, Lemmas \ref{med} and \ref{v2} imply that there is
     either a root $\alpha$ of $F(Z,1)$ with \eqref{med2} or a root $\alpha^{-1}$ of $F(1,Z)$ with \eqref{v22}.
     Further, the right hand sides of these inequalities increase with $u$ and $v$ respectively. Therefore
     we may replace $u,v$ with $s$. This, together with Lemma \ref{set_S}, gives the assertion of the lemma.
    \end{proof}

\section{Estimation of medium solutions}
Let $Y_W$ be as given in \eqref{YL}. In this section, we shall estimate $P_{med}(Y_W,Y)$ for suitably chosen $Y$.     
 \begin{lemma}\label{m1}
Let $F(X,Y)$ be given by \eqref{poly}.
\begin{enumerate}
 \item[(i)] Let
  \[
       Y_S=((12e^{\Psi})^{r}R_1^{2s}h)^{\frac{1}{r-2s}}.
      \]
Then
      \[
     P_{med}(Y_W,Y_S)  \ll \frac{s}{\Phi}(\log s+\log(1+\log h^{1/r})).
      \]
 \item[(ii)]  Suppose that the coefficients of $F(X,Y)$ satisfy \eqref{st_line_1} and that $r\geq 4s$. Let
\[
 Y_S'=(8^{r}R_1^{s}(s^2r)^{3s}h)^{\frac{1}{r-2s}}.
\]
Then
\[
 P_{med}(Y_W,Y_S')\ll s(\log s+\log(1+\log h^{1/r})).
\]
      \end{enumerate}
 \end{lemma}
 \begin{proof}\noindent 
  $(i)$ By Lemma \ref{app_lemma}, it is enough to estimate the number of primitive pairs $(x,y)$ satisfying \eqref{medium_app}
  for some $\alpha\in S_1$ or \eqref{29dash} for some $\alpha^*\in S_1^*$
  with
      \[
       Y_S\leq y\leq Y_W.
      \]
We consider the case when \eqref{medium_app} is satisfied. The other case is similar. Let $\alpha\in S_1$ and
let $U=\{(x_1,y_1),\ldots,(x_{\nu},y_{\nu})\}$ be the set of all solutions of \eqref{medium_app} with $\gcd(x_i,y_i)=1$ and
      \[
       Y_S\leq y_1\leq\ldots\leq y_{\nu}\leq Y_W.
      \]
Suppose that $\nu\geq 2$. Then
      \begin{eqnarray*}
         \frac{1}{y_iy_{i+1}}&\leq&\left|\frac{x_i}{y_i}-\frac{x_{i+1}}{y_{i+1}}\right|\leq 
         \left|\alpha-\frac{x_i}{y_i}\right|+\left|\alpha-\frac{x_{i+1}}{y_{i+1}}\right|\\
         &\leq& \frac{K_1}{2y_i^{r/s}}+\frac{K_1}{2y_{i+1}^{r/s}}\leq\frac{K_1}{y_i^{r/s}},
      \end{eqnarray*}
where 
      \[
       K_1=2R_1(rs)^2(12e^{\Psi})^{r/s}h^{1/s}H^{(1/r)-(1/s)}.
      \]
Thus we have
      \begin{eqnarray*}
       y_{i+1}&\geq& K_1^{-1}y_i^{(r/s)-1}.
      \end{eqnarray*}
In Lemma \ref{count_lemma}(i) take $T((x_i,y_i))=y_i$, $\beta=\frac{1}{K_1}$, $\gamma=\frac{r-s}{s}$, $A_1=Y_S$
and $B_1=Y_W$. Note that $\gamma=(r/s)-1>4e^{2\Phi}-1\geq e^{\Phi}\geq 2$. Further, since $R_1\geq 4(rs)^4$, we have
\[
 \frac{Y_S}{K_1^{\frac{1}{\kappa(\gamma-1)}}}\geq \frac{R_1^{\frac{s}{r-2s}}H^{\frac{r-s}{2r(r-2s)}}}{(rs)^{\frac{2s}{r-2s}}2^{\frac{s}{r-2s}}}
 \geq R_1^{\frac{s}{2(r-2s)}}H^{\frac{r-s}{2r(r-2s)}}.
\]
Using the inequality $M\leq (r+1)H$ (\cite[Eqn $(6)$]{Ma0}), we get
\[
 \log Y_W\ll \sqrt{r}+\log H+\log h^{1/r}.
\]
Now we apply Lemma \ref{count_lemma}(i) to obtain 
     \begin{eqnarray*}
      \nu&\ll&1+\frac{1}{\log\gamma}\log\left(\frac{2r(r-2s)(\sqrt{r}+\log H+\log h^{1/r})}{(r-s)\log H+rs\log R_1}\right)\\
      &\ll&\frac{\log r+\log(1+\log h^{1/r})}{\log\gamma}.
     \end{eqnarray*}
 If $r\ll s^3$, we get
       \begin{equation*}\label{nu_temp}
        \nu\ll\frac{\log s+\log(1+\log h^{1/r})}{\Phi}.
       \end{equation*}
If $r\gg s^3$, then $\gamma\gg r^{2/3}$, which implies that
      \[
         \nu\ll 1+\frac{\log(1+\log h^{1/r})}{\log r}\ll\frac{\log s+\log(1+\log h^{1/r})}{\Phi}
      \]
as $\Phi\ll\log s$ (see Remark \ref{rk_13}).
Since this is true for each $\alpha\in S_1$, we obtain the assertion of part $(i)$ of the lemma.\\
$(ii)$ By Remark \ref{rk_lemma} and Lemmas \ref{dvt1}, \ref{existence_u_v} and \ref{set_S}, there is 
a set $S_1\subseteq S$ with $|S_1|\ll s$,
 such that for some $\alpha\in S_1$,
    \begin{equation}\label{d_temp}
      d\left(S_1,\frac{x}{y}\right)=\left|\alpha-\frac{x}{y}\right|\leq\frac{rR_1}{2}\left(\frac{s(rs^2)^{s-1}2^{r+s+1}h}{|y|^r|a_s||\alpha|^{r-u}}\right)^{1/u}.
    \end{equation}
By Lemma \ref{bds}, it follows that every root $\alpha$ of $f$ satisfies
    \[
     \frac{1}{2}e^{\sigma}<|\alpha|<2e^{\sigma},
    \]
where $\sigma$
is the slope of the line joining $(0,-\log |a_0|)$ and $(r,-\log |a_s|)$. This implies that
\[
 \log |\alpha|>\sigma-\log 2=\frac{-\log |a_s|+\log |a_0|}{r}-\log 2.
\]
Therefore
\[
 |a_s||\alpha|^{r-u}\geq \frac{|a_s|^{u/r}|a_0|^{(r-u)/r}}{2^{r-u}}.
\]
By \eqref{d_temp} and $|y|\geq Y_S'$, we get that
\[
       d\left(S_1,\frac{x}{y}\right)\leq\frac{rR_1}{2}\left(\frac{s(rs^2)^{s-1}2^{3r}h}{|y|^r|a_s|^{s/r}|a_0|^{(r-s)/r}}\right)^{1/s}.
\]
As the height of $F$ is either $|a_0|$
     or $|a_s|$, we obtain that
\begin{equation}\label{35}
       \left|\alpha-\frac{x}{y}\right|\leq\frac{rR_1}{2H^{1/r}}\left(\frac{s(rs^2)^{s-1}2^{3r}h}{|y|^r}\right)^{1/s}.
\end{equation}
Let $U=\{(x_1,y_1),\ldots,(x_{\nu},y_{\nu})\}$ be the set of all the solutions of \eqref{35} with $\gcd(x_i,y_i)=1$ and
      \[
       Y_S'\leq y_1\leq\ldots\leq y_{\nu}\leq Y_W.
      \]
Put
      \[
        K_2=\frac{R_18^{r/s}(sh)^{1/s}(rs)^2}{H^{1/r}}.
       \]
Suppose that $\nu\geq 2$.
As in the proof of $(i)$, we get from \eqref{35} that
\[
 y_{i+1}\geq K_2^{-1}y_i^{(r/s)-1}.
\]
In Lemma \ref{count_lemma}(i) take $T((x_i,y_i))=y_i$, $\beta=\frac{1}{K_2}$, $\gamma=\frac{r-s}{s}$, $A_1=Y_S'$
and $B_1=Y_W$. Note that
  \[
        Y_S'\geq K_2^{\frac{1}{\kappa(\gamma-1)}}H^{\frac{s}{2r(r-2s)}}(rs^3)^{\frac{s}{r-2s}} \textrm{ and } \gamma\geq 2.
       \]
Thus by Lemma \ref{count_lemma}(i), we obtain
 \begin{eqnarray*}
      \nu&\ll&1+\frac{1}{\log\gamma}\log\left(\frac{2r(r-2s)(\sqrt{r}+\log H+\log h^{1/r})}{s\log H+2rs\log(rs^3)}\right)\\
      &\ll&\frac{\log r+\log(1+\log h^{1/r})}{\log\gamma}.
     \end{eqnarray*}
Now we argue as in the proof of $(i)$ to get the assertion.
\end{proof}
\section{Estimation of small solutions}
To estimate the number of small solutions, we use the following lemma from \cite{MS2}.
\begin{lemma}\cite[Lemma 18]{MS2}\label{small}\\
Let $F(X,Y)$ be given by \eqref{poly} and let $r\geq 4s$.
Then for any $Y\geq 1$, we have
      \[
    P_{sma}(Y)  \ll (rs^2)^{2s/r} h^{2/r} + s\ Y.
      \]
\end{lemma}
\begin{lemma}\label{8.3}
 Let $F(X,Y)$ be given by \eqref{poly}.
 Then
 \begin{enumerate}
  \item[(i)] $$ P_{sma}(Y_S) \ll e^{c_7(\log s) e^{-2\Phi}}h^{\frac{2}{r}}+se^{\Phi+c_8(\log^3s) e^{-\Phi}}h^{\frac{1}{r-2s}}.$$
  \item[(ii)]
 \[
  P_{sma}(Y_S')\ll sh^{2/r} \textrm{ whenever } r\geq s\log^3s.
 \]
 \end{enumerate}
\end{lemma}
\begin{proof}
 $(i)$ follows from \eqref{assumption_r} and Lemma \ref{small} with $Y$ as $Y_S$. Similarly, $(ii)$ follows by taking $Y$ as $Y_S'$ in Lemma \ref{small}
 and using $r\geq s\log^3s$.
\end{proof}
\section{Proofs of Proposition \ref{prop} and Theorems \ref{thm_1}, \ref{line}}\label{last}
\noindent \textbf{Proof of Proposition \ref{prop}}\\
Suppose that $r\geq 4s$. Let $P(h)$ denote the number of primitive solutions of \eqref{Thue_ineq}. Then
\[
 P(h)=P_{\ell ar}(Y_W)+P_{med}(Y_W,Y_S)+P_{sma}(Y_S).
\]
The upper bounds for the three quantities on the right hand side are obtained from \eqref{la}, Lemma \ref{m1}(i)
and Lemma \ref{8.3}(i). This together with $\Phi\ll\log s$ (see Remark \ref{rk_13}) yields
\[
 P(h)\ll \rho h^{2/r},
\]
where 
\[
\rho=\frac{s \log s}{\Phi}+s^{c_7e^{-2\Phi}}+se^{\Phi+c_8(\log^3s) e^{-\Phi}}.
\]
Using a partial summation argument, it was shown in \cite[p. $212$]{MS2} that
\begin{equation}\label{part_sum}
\nonumber N_F(h)\ll P(h)+h^{1/r}r^{-1}\sum_{n=1}^{h-1} P(n)n^{-1-(1/r)}.
\end{equation}
Substituting our estimate for $P(h)$, we obtain that
\[
 N_F(h)\ll \rho h^{2/r}.
\]
When $r<4s$, we use \eqref{Schmidt_result} to obtain 
\[
 N_F(h)\ll s\ C_1(r,h).
\]
This proves the proposition.
\\
\textbf{Proof of Theorem \ref{thm_1}}\\
Observe that $\Phi\geq 3\log\log s$. Thus from Proposition \ref{prop}, we get
\[
                N_F(h)\ll s(\log s+e^{\Phi}) C_1(r,h)\ll se^{\Phi} C_1(r,h).
\]
\\
\textbf{Proof of Theorem \ref{line}}\\
We have
\[
 P(h)=P_{\ell ar}(Y_W)+P_{med}(Y_W,Y_S')+P_{sma}(Y_S').
\]
We use the respective upper bounds for the three quantities on the right hand side from \eqref{la}, 
Lemma \ref{m1}(ii) and Lemma \ref{8.3}(ii) and proceed as in the proof of Proposition \ref{prop} to give the assertion.
 
\end{document}